\documentclass{amsart}
\usepackage{amssymb}
\usepackage{amsmath}
\newtheorem{theorem}{Theorem}[section]

\newcommand{\mpn}{\medskip\par\noindent}
\newcommand{\pn}{\par\noindent}

\newcommand{\bpn}{\bigskip\par\noindent}
\theoremstyle{definition}

\numberwithin{equation}{section}

\begin{document}
\newcommand{\Mod}[1]{\,(\text{\mbox{\rm mod}}\;#1)}
\title[  A note on fermionic $p$-adic integrals on $\Bbb Z_p$ and Umbral calculus ]{ A note on fermionic $p$-adic integrals on $\Bbb Z_p$ and Umbral calculus}
\author{Taekyun Kim, Dae San Kim, Sangtae Jeong and Seog-Hoon Rim}
\begin{abstract}
In this paper  we study some properties of the fermionic $p$-adic integrals on $\Bbb Z_p$ arising from the umbral calculus.
\end{abstract}
\maketitle
\vspace{3mm}

\def\ord{\text{ord}_p}
\par\bigskip\noindent
\section{Introduction }
 Let $p$ be a fixed odd prime number. Throughout this paper
$\mathbb{Z}_{p},\mathbb{Q}_{p}$  and $\mathbb{C}_{p}$
denote the ring of $p$-adic rational integers, the field of $p$-adic
rational numbers  and the completion of
the algebraic closure of $\mathbb{Q}_{p}$, respectively. Let
$\mathbb{N}$ be the set of natural numbers and
$\mathbb{Z}_{+}=\mathbb{N}\cup\{0\}$. Let $C(\mathbb{Z}_{p})$ be the space of continuous functions on $\mathbb{Z}_{p}$.
For $f \in C(\mathbb{Z}_{p}) $, the fermionic $p$-adic integral on $\Bbb Z_p$ is defined by
\begin{equation}\tag{1}
\begin{split}
 \int_{\mathbb{Z}_{p}}f(x)d\mu_{-1}(x)&=\lim_{N
\rightarrow \infty }  \sum_{x=0}^{p^{N}-1}
f(x)\mu_{-1}(x+p^N\mathbb{Z}_{p})
\\&=\lim_{N
\rightarrow \infty }  \sum_{x=0}^{p^{N}-1}
f(x)(-1)^x, \text{ (see [1,2,11])}.
\end{split}
\end{equation}
For $n \in \mathbb{N}$, we have
\begin{equation}\tag{2}
\begin{split}
 \int_{\mathbb{Z}_{p}}f(x+n)d\mu_{-1}(x)+(-1)^{n-1} \int_{\mathbb{Z}_{p}}f(x )d\mu_{-1}(x)=2\sum_{l=0}^{n-1}(-1)^{n-1-l}f(l).
\end{split}
\end{equation}
In the special case, $n=1$, we note that
\begin{equation}\tag{3}
\begin{split}
 \int_{\mathbb{Z}_{p}}f(x+1)d\mu_{-1}(x)+  \int_{\mathbb{Z}_{p}}f(x )d\mu_{-1}(x)=2 f(0),  \text{ (see [ 11])}.
\end{split}
\end{equation}
Let $\mathcal{F}$ be the set of all formal power series in the variable $t$ over $\mathbb{C}_{p}$ with
\begin{equation*}
\begin{split}
\mathcal{F}=\Bigl{\{} f(t)=\sum_{k=0}^{\infty} \frac{a_k}{k!}t^k \Bigl{|} a_k \in \mathbb{C}_{p} \Bigl{\}}.
\end{split}
\end{equation*}
Let $\mathbb{ P}=\mathbb{C}_{p}[x]$ and let $\mathbb{ P}^*$ denote the vector space of all linear functionals on $\mathbb{ P}$.
The formal power series
\begin{equation}\tag{4}
\begin{split}
f(t)=\sum_{k=0}^{\infty} \frac{a_k}{k!}t^k \in \mathcal{F}.
\end{split}
\end{equation}
defines a linear functional on   $\mathbb{ P}$ by setting
\begin{equation}\tag{5}
\begin{split}
\langle f(t) | x^n \rangle=a_n  \text{ for all } n\geq 0,  \text{ (see [ 7,14])}.
\end{split}
\end{equation}
Thus, by (4) and (5), we get
\begin{equation}\tag{6}
\begin{split}
\langle t^k | x^n \rangle=n!\delta_{n,k},~(n,~k\geq 0),
\end{split}
\end{equation}
where $\delta_{n,k}$ is the Kronecker symbol (see [7,14]).
Here, $\mathcal{F}$ denotes both the algebra of formal power series in $t$ and the vector space of all linear functionals on  $\mathbb{ P}$,
and so an element $f(t)$ of $\mathcal{F}$ will be thought of as both a formal power series and a linear functional. We shall call $\mathcal{F}$
the umbral algebra. The umbral calculus is the study of umbral algebra (see [7,14]).
\\
The order $O(f(t))$ of power series $f(t)(\neq0)$ is the smallest integer $k$ for which $a_k$ does not vanish (see [7,4]). The series $f(t)$
has a multiplicative inverse, denoted by $f(t)^{-1}$ or $\frac{1}{f(t)}$, if and only if $O(f(t))=0$.
Such series is called an invertible series. A series $f(t)$ for which $O(f(t))=1$ is called a delta series (see [7,14]). For $f(t),~g(t) \in \mathcal{F}$,
we have $\langle f(t)g(t)|p(x)\rangle=\langle f(t)|g(t)p(x)\rangle=\langle g(t)|f(t)p(x)\rangle$. By (6), we get
\begin{equation}\tag{7}
\begin{split}
\langle e^{yt} | x^n \rangle=y^n, \quad \langle e^{yt} |p(x)\rangle=p(y),  \text{ (see [ 7,14])}.
\end{split}
\end{equation}
Let $f(t) \in \mathcal{F}$. Then we note that
\begin{equation}\tag{8}
\begin{split}
f(t)=\sum_{k=0}^{\infty}\frac{\langle f(t) | x^k \rangle}{k!}t^k,
\end{split}
\end{equation}
and
\begin{equation}\tag{9}
\begin{split}
p(x)=\sum_{k=0}^{\infty}\frac{\langle t^k | p(x) \rangle}{k!}x^k,\text{ for } p(x) \in  \mathbb{P}, \text{ (see [14])}.
\end{split}
\end{equation}
Let $f_{1}(t),~f_{2}(t),~\cdots   ,~f_{m}(t) \in \mathcal{F}$. It is known in [7,14] that
\begin{equation}\tag{10}
\begin{split}
\langle f_{1}(t)\cdots f_{m}(t)   | x^n \rangle =\sum \binom{n}{i_{1},\cdots,i_{m} } \langle f_{1}(t) | x^{i_{1}}\rangle \cdots \langle f_{m}(t) | x^{i_{m}} \rangle,
\end{split}
\end{equation}
where the sum is over all nonnegative integers $i_{1},~\cdots,~i_{m}$ such that $i_{1}+i_{2}+\cdots+i_{m}=n$ (see [7,14]).

By (9), we get
\begin{equation}\tag{11}
\begin{split}
p^{(k)}(x)=\frac{d^kp(x)}{dx^k}&=\sum_{l=k}^{\infty}\frac{\langle t^l | p(x) \rangle}{l!}l(l-1)\cdots(l-k+1)x^{l-k}
\\&=\sum_{l=k}^{\infty}\langle t^l | p(x) \rangle \binom{l}{k} \frac{k!}{l!}x^{l-k}.
\end{split}
\end{equation}
Thus, from (11), we have
\begin{equation}\tag{12}
\begin{split}
p^{(k)}(0)=\langle t^k | p(x) \rangle =\langle 1 | p^{(k)}(x) \rangle ,
\end{split}
\end{equation}
and
\begin{equation}\tag{13}
\begin{split}
t^k   p(x) = p^{(k)}(x)=\frac{d^kp(x)}{dx^k} \quad \text{ (see [7,14])}.
\end{split}
\end{equation}
From (13), we note that
\begin{equation}\tag{14}
\begin{split}
e^{yt}p(x)=p(x+y)\quad \text{ (see [7,14])}.
\end{split}
\end{equation}
In this paper, $s_n(x)$ denotes a polynomial of degree $n$. Let us assume that $f(t),~g(t) \in \mathcal{F}$
with $o(f(t))=1$ and $o(g(t))=1$. Then there exists a unique sequence $s_n(x)$ of polynomials satisfying
$\langle g(t)f(t)^k |s_n(x)\rangle=n!\delta_{n,k}$ for all $n,~k\geq 0$.
The sequence $s_n(x)$ is called the Sheffer sequence for $(g(t),~f(t))$, which is denoted by $s_n(x)\sim (g(t),~f(t))$.
If  $s_n(x)\sim (g(t),~t)$, then $s_n(x)$ is called the Appell sequence for $g(t)$ (see [7,14]).

Let $p(x) \in \mathbb{P}$. Then we note that
\begin{equation}\tag{15}
\begin{split}
 \langle f(t) | xp(x) \rangle = \langle  \partial _tf(t) |  p(x) \rangle  = \langle  f^{\prime}(t) |  p(x) \rangle,
\end{split}
\end{equation}
and
\begin{equation}\tag{16}
\begin{split}
 \langle e^{yt}-1 | p(x) \rangle = p(y)-p(0), \quad \text{ (see [7,14])}.
\end{split}
\end{equation}
Let us assume that $s_n(x)\sim (g(t),~f(t))$. Then we have
\begin{equation}\tag{17}
\begin{split}
h(t)=\sum_{k=0}^{\infty}\frac{\langle h(t) |s_k(x) \rangle}{k!}g(t)f(t)^k,\quad h(t) \in \mathcal{F},
\end{split}
\end{equation}
\begin{equation}\tag{18}
\begin{split}
p(x)=\sum_{k=0}^{\infty}\frac{\langle g(t)f(t)^k |p(x) \rangle}{k!}s_k(x),\quad p(x) \in \mathbb{P},
\end{split}
\end{equation}
\begin{equation}\tag{19}
\begin{split}
\frac{1}{g(\bar{f}(t))}e^{y\bar{f}(t)}=\sum_{k=0}^{\infty}\frac{s_k(y)}{k!}t^k,\quad \text{for all } y \in \mathbb{C}_p,
\end{split}
\end{equation}
where $\bar{f}(t)$ is the compositional inverse of $f(t)$, and
\begin{equation}\tag{20}
\begin{split}
f(t)s_n(x)=ns_{n-1}(x), \quad \text{ (see [7,14])}.
\end{split}
\end{equation}
As is well known, the Euler polynomials are defined by the generating function to be
\begin{equation}\tag{21}
\begin{split}
\frac{2}{e^t+1}e^{xt}=e^{E(x)t}=\sum_{n=0}^{\infty}E_n(x)\frac{t^n}{n!}, \quad  \text{ (see [1-19])},
\end{split}
\end{equation}
with the usual convention about replacing $E^n(x)$ by $E_n(x)$. In the special case, $x=0$, $E_n(0)=E_n $ are called the $n$-th Euler numbers

Let $s_n(x)\sim (g(t),~t)$. Then Appell identity is known to be
\begin{equation}\tag{22}
\begin{split}
s_n(x+y)=\sum_{k=0}^{n}\binom{n}{k}s_{n-k}(x)y^{k}=\sum_{k=0}^{n}\binom{n}{k}s_{ k}(x)y^{n-k}.
\end{split}
\end{equation}

From (21), we note that the recurrence relation of the Euler numbers is given by
\begin{equation}\tag{23}
\begin{split}
E_0=1, \quad (E+1)^n+E_n=E_n(1)+E_n=2\delta_{0,n}.
\end{split}
\end{equation}
By (1) and (21), we get
\begin{equation}\tag{24}
\begin{split}
\int_{\mathbb{Z}_p}(x+y)^n d\mu_{-1}(y)=E_n(x), \quad \int_{\mathbb{Z}_p}x^n d\mu_{-1}(y)=E_n,
\end{split}
\end{equation}
where $n\geq 0$ (see [1,11,16]).

Recently, R. Dere and Y. Simsek have studied applications of umbral algebra to some special polynomials (see [7]).
In this paper we study some properties of the fermionic $p$-adic integrals on $\mathbb{Z}_p$ arising from the umbral calculus.
\par\bigskip\noindent
\section{ Umbral calculus and fermionic $p$-adic integrals on $\mathbb{Z}_p$ }

Let $s_n(x)\sim (g(t),~t)$. Then, by (19), we get
\begin{equation}\tag{25}
\begin{split}
\frac{1}{g(t)}x^n= s_{n }(x) \quad \text{if and only if} \quad x^n=g(t)s_{n }(x) .
\end{split}
\end{equation}
Let us assume that $g(t)=\frac{e^t+1}{2}$. Then we note that $g(t)$ is an invertible functional.
\\
By (21), we get
\begin{equation}\tag{26}
\begin{split}
\frac{1}{g(t)}e^{xt} =\sum_{k=0}^{\infty} E_k(x) \frac{t^k}{k!}.
\end{split}
\end{equation}
Thus, from (26), we have
\begin{equation}\tag{27}
\begin{split}
\frac{1}{g(t)}x^n= E_{n }(x), \quad  tE_{n }(x)=\frac{n}{g(t)}x^{n-1}=n E_{n -1}(x).
\end{split}
\end{equation}
By (19), (20) and (27), we see that $E_{n }(x)$ is an Appell sequence for $g(t)=\frac{e^t+1}{2}$.

It is easy to show that
\begin{equation}\tag{28}
\begin{split}
  E_{n +1}(x)=\Bigl(x-\frac{g^{\prime}(t)}{g(t)}\Bigl)E_{n }(x), \quad (n\geq0).
\end{split}
\end{equation}
From (2), (21) and (24), we note that
\begin{equation}\tag{29}
\begin{split}
\int_{\mathbb{Z}_p}e^{(x+y+1)t}d\mu_{-1}(y)+ \int_{\mathbb{Z}_p}e^{(x+y )t} d\mu_{-1}(y)=2e^{xt}.
\end{split}
\end{equation}
Thus, by (29), we get
\begin{equation}\tag{30}
\begin{split}
\int_{\mathbb{Z}_p}(x+y+1)^nd\mu_{-1}(y)+ \int_{\mathbb{Z}_p}(x+y )^nd\mu_{-1}(y)=2x^n.
\end{split}
\end{equation}
From (24) and (30), we have
\begin{equation}\tag{31}
\begin{split}
  E_{n }(x+1)+ E_{n }(x )=2x^n, \quad (n\geq0).
\end{split}
\end{equation}
By (28), we see that
\begin{equation}\tag{32}
\begin{split}
 g(t) E_{n+1 }(x )=g(t)xE_{n }(x )-g^{\prime}(t)E_{n }(x ), \quad (n\geq0).
\end{split}
\end{equation}
Thus, we have
\begin{equation}\tag{33}
\begin{split}
(e^t+1) E_{n+1 }(x )=(e^t+1)xE_{n }(x )-e^tE_{n }(x ).
\end{split}
\end{equation}
By (33), we get
\begin{equation}\tag{34}
\begin{split}
  E_{n+1 }(x+1 )+ E_{n+1 }(x )=(x+1) E_{n }(x+1 )+xE_{n }(x )-E_{n }(x+1 ).
\end{split}
\end{equation}
Thus, from (34) and (31), we have
\begin{equation}\tag{35}
\begin{split}
  E_{n+1 }(x+1 )+ E_{n+1 }(x )=x( E_{n }(x+1 )+ E_{n }(x)).
\end{split}
\end{equation}
By (35), we get
\begin{equation*}
\begin{split}
  E_{n }(x+1 )+ E_{n  }(x )&=x( E_{n-1 }(x+1 )+ E_{n-1 }(x))=x^2(E_{n-2 }(x+1 )+ E_{n-2 }(x))
  \\&=\cdots=x^n(E_{0 }(x+1 )+ E_{0 }(x))=2x^n.
\end{split}
\end{equation*}
Let us consider the functional $f(t)$ such that
\begin{equation}\tag{36}
\begin{split}
\langle f(t)|p(x) \rangle = \int_{\mathbb{Z}_p}p(u)d\mu_{-1}(u),
\end{split}
\end{equation}
for all polynomials $p(x)$. It can be determined from (8) to be
\begin{equation}\tag{37}
\begin{split}
f(t) =\sum_{k=0}^{\infty} \frac{\langle f(t)|x^k \rangle}{k!} t^k = \sum_{k=0}^{\infty}\int_{\mathbb{Z}_p}u^k d\mu_{-1}(u) \frac{t^k}{k!}
 =\int_{\mathbb{Z}_p}e^{ut} d\mu_{-1}(u).
\end{split}
\end{equation}
By (29) and (37), we get
\begin{equation}\tag{38}
\begin{split}
f(t)  =\int_{\mathbb{Z}_p}e^{ut} d\mu_{-1}(u)=\frac{2}{e^t+1}.
\end{split}
\end{equation}
Therefore, by (38), we obtain the following theorem.
\begin{theorem}
For $p(x) \in \mathbb{P}$, we have
\begin{equation*}
\langle \int_{\mathbb{Z}_p}e^{yt} d\mu_{-1}(y) | p(x) \rangle = \int_{\mathbb{Z}_p}p(u) d\mu_{-1}(u).
\end{equation*}
That is,
\begin{equation*}
\langle \frac{2}{e^t+1} | p(x) \rangle = \int_{\mathbb{Z}_p}p(u) d\mu_{-1}(u).
\end{equation*}
Also, the $n$-th Euler number is given by
\begin{equation*}
E_n=\langle\int_{\mathbb{Z}_p}e^{yt} d\mu_{-1}(y) | x^n \rangle .
\end{equation*}
\end{theorem}
\vspace{3mm}
By (3) and (30), we get
\begin{equation}\tag{39}
\begin{split}
\sum_{n=0}^{\infty}\int_{\mathbb{Z}_p}(x+y)^n d\mu_{-1}(y)\frac{t^n}{n!}=\int_{\mathbb{Z}_p}e^{(x+y)t} d\mu_{-1}(y)
=\sum_{n=0}^{\infty} \int_{\mathbb{Z}_p}e^{yt} d\mu_{-1}(y)x^n\frac{t^n}{n!}.
\end{split}
\end{equation}
From (24) and (39), we have
\begin{equation}\tag{40}
\begin{split}
E_n(x)=\int_{\mathbb{Z}_p}e^{yt} d\mu_{-1}(y)x^n=\frac{2}{e^t+1}x^n,
\end{split}
\end{equation}
where $n\geq0$.

Therefore, by (40), we obtain the following theorem.
\begin{theorem}
For $p(x) \in \mathbb{P}$, we have
\begin{equation*}
\int_{\mathbb{Z}_p}p(x+y) d\mu_{-1}(y)=\int_{\mathbb{Z}_p}e^{yt} d\mu_{-1}(y)p(x)=\frac{2}{e^t+1}p(x).
\end{equation*}
\end{theorem}

From (22), we note that
\begin{equation*}
E_n(x+y)=\sum_{k=0}^{n} \binom{n}{k}E_k(x)y^{n-k}.
\end{equation*}
The Euler polynomials of order $r$ are defined by the generating function to be
\begin{equation}\tag{41}
\begin{split}
\underbrace{\Bigl(\frac{2}{e^t+1}\Bigl)\times\Bigl(\frac{2}{e^t+1}\Bigl)\times \cdots \times \Bigl(\frac{2}{e^t+1}\Bigl) }_{r\text{
times}} e^{xt}=\Bigl(\frac{2}{e^t+1}\Bigl)^re^{xt}
=\sum_{n=0}^{\infty}E_n^{(r)}(x)\frac{t^n}{n!}.
\end{split}
\end{equation}
In the special case, $x=0$,
$E_n^{(r)}(0)=E_n^{(r)}$ are called the $n$-th Euler numbers of order $r$ $(r\geq0)$, (see [1-19]). Let us take $g^r(t)=\Bigl(\frac{e^t+1}{2}\Bigl)^r$.
Then we see that $g^r(t)$ is an invertible functional in $\mathcal{F}$. By (41), we get
\begin{equation}\tag{42}
\begin{split}
\frac{1}{g^r(t)}e^{xt}
=\sum_{n=0}^{\infty}E_n^{(r)}(x)\frac{t^n}{n!}.
\end{split}
\end{equation}
Thus, we have
\begin{equation}\tag{43}
\begin{split}
\frac{1}{g^r(t)}x^n
=E_n^{(r)}(x), \quad tE_n^{(r)}(x)=\frac{n}{g^r(t)}x^{n-1}=nE_{n-1}^{(r)}(x).
\end{split}
\end{equation}
So, by (42), we see that $E_n^{(r)}(x)$ is the Appell sequence for $ \Bigl(\frac{e^t+1}{2}\Bigl)^r$. From (22), we have
\begin{equation}\tag{44}
\begin{split}
 E_n^{(r)}(x+y)=\sum_{k=0}^{n} \binom{n}{k}E_{n-k}^{(r)}(x)y^{ k}.
\end{split}
\end{equation}
It is easy to show that
\begin{equation}\tag{45}
\begin{split}
\int_{\mathbb{Z}_p}\cdots\int_{\mathbb{Z}_p} e^{(x_1+\cdots x_r+x)t} d\mu_{-1}(x_1)\cdots d\mu_{-1}(x_r)=\Bigl(\frac{2}{e^t+1}\Bigl)^re^{xt}
=\sum_{n=0}^{\infty}E_n^{(r)}(x)\frac{t^n}{n!}.
\end{split}
\end{equation}
By (6) and (45), we get
\begin{equation}\tag{46}
\begin{split}
E_n^{(r)}=  \langle \underbrace{\int_{\mathbb{Z}_p}\cdots\int_{\mathbb{Z}_p}}_{r-\text{times}} e^{(x_1+\cdots x_r )t} d\mu_{-1}(x_1)\cdots d\mu_{-1}(x_r) |x^n \rangle ,
\quad (n\geq0),
\end{split}
\end{equation}
and, by (10),
\begin{equation}\tag{47}
\begin{split}
  &\langle \int_{\mathbb{Z}_p}\cdots\int_{\mathbb{Z}_p} e^{(x_1+\cdots x_r )t} d\mu_{-1}(x_1)\cdots d\mu_{-1}(x_r) |x^n \rangle
  \\&=\sum_{n=i_1+\cdots + i_r} \binom {n}{i_1,\cdots ,i_r}\langle \int_{\mathbb{Z}_p}  e^{ x_1 t} d\mu_{-1}(x_1)  |x^{i_1} \rangle \cdots \langle \int_{\mathbb{Z}_p}  e^{  x_r  t} d\mu_{-1}(x_r)  |x^{i_r} \rangle
  \\&=\sum_{n=i_1+\cdots + i_r} \binom {n}{i_1,\cdots ,i_r}E_{i_1}E_{i_2}\cdots E_{i_r}.
\end{split}
\end{equation}
From (46) and (47), we have
\begin{equation}\tag{48}
\begin{split}
E_n^{(r)}= \sum_{n=i_1+\cdots + i_r} \binom {n}{i_1,\cdots ,i_r}E_{i_1} \cdots E_{i_r}.
\end{split}
\end{equation}
By (44) and (48), we see that $E_n^{(r)}(x)$ is  a monic polynomial of degree $n$ with coefficients in $\mathbb{Q}$.
Let $r \in \mathbb{N}$. Then we note that
\begin{equation}\tag{49}
\begin{split}
g^r(t)= \frac{1}{\underbrace{\int_{\mathbb{Z}_p}\cdots\int_{\mathbb{Z}_p}}_{r-\text{times}} e^{(x_1+\cdots x_r )t} d\mu_{-1}(x_1)\cdots d\mu_{-1}(x_r) }=\Bigl(\frac{e^t+1}{2}\Bigl)^r.
\end{split}
\end{equation}
By (49), we get
\begin{equation}\tag{50}
\begin{split}
\frac{1}{g^r(t)}e^{xt}=  \int_{\mathbb{Z}_p}\cdots\int_{\mathbb{Z}_p} e^{(x_1+\cdots x_r +x)t} d\mu_{-1}(x_1)\cdots d\mu_{-1}(x_r)  =\sum_{n=0}^{\infty}E_n^{(r)}(x)\frac{t^n}{n!}.
\end{split}
\end{equation}
From (50), we have
\begin{equation}\tag{51}
\begin{split}
E_n^{(r)}(x)&=  \int_{\mathbb{Z}_p}\cdots\int_{\mathbb{Z}_p}  (x_1+\cdots x_r +x)^n d\mu_{-1}(x_1)\cdots d\mu_{-1}(x_r)
\\&=\frac{1}{\underbrace{\int_{\mathbb{Z}_p}\cdots\int_{\mathbb{Z}_p}}_{r-\text{times}} e^{(x_1+\cdots x_r )t} d\mu_{-1}(x_1)\cdots d\mu_{-1}(x_r) }x^n
\\&=\frac{1}{g^r(t)}x^n.
\end{split}
\end{equation}
Therefore, by (51), we obtain the following theorem.
\begin{theorem}
For $p(x) \in \mathbb{P}$ and $r \in \mathbb{N}$. Then we have
\begin{equation*}
 \underbrace{\int_{\mathbb{Z}_p}\cdots\int_{\mathbb{Z}_p}}_{r-\text{times}} p (x_1+\cdots x_r +x)  d\mu_{-1}(x_1)\cdots d\mu_{-1}(x_r)=\Bigl(\frac{2}{e^t+1}\Bigl)^rp(x).
\end{equation*}
In particular,
\begin{equation*}
 E_n^{(r)}(x)=\Bigl(\frac{2}{e^t+1}\Bigl)^rx^n=\int_{\mathbb{Z}_p}\cdots\int_{\mathbb{Z}_p} e^{(x_1+\cdots x_r )t} d\mu_{-1}(x_1)\cdots d\mu_{-1}(x_r)x^n.
\end{equation*}
That is,
\begin{equation*}
 E_n^{(r)}(x)\sim \Bigl(\frac{1}{\int_{\mathbb{Z}_p}\cdots\int_{\mathbb{Z}_p} e^{(x_1+\cdots x_r )t} d\mu_{-1}(x_1)\cdots d\mu_{-1}(x_r) },t\Bigl).
\end{equation*}
\end{theorem}
Let us take the functional $f^r(t)$ such that
\begin{equation}\tag{52}
\begin{split}
\langle f^r(t) | p(x) \rangle =\underbrace{\int_{\mathbb{Z}_p}\cdots\int_{\mathbb{Z}_p}}_{r-\text{times}} p (x_1+\cdots x_r  )  d\mu_{-1}(x_1)\cdots d\mu_{-1}(x_r),
\end{split}
\end{equation}
for all polynomials $p(x)$. It can be determined from (8) to be
\begin{equation}\tag{53}
\begin{split}
 f^r(t)&=\sum_{k=0}^{\infty} \frac{\langle f^r(t) |x^k\rangle}{k!}t^k
 \\&=\sum_{k=0}^{\infty} \underbrace{\int_{\mathbb{Z}_p}\cdots\int_{\mathbb{Z}_p}}_{r-\text{times}}  (x_1+\cdots x_r  )^k  d\mu_{-1}(x_1)\cdots d\mu_{-1}(x_r)
 \frac{t^k}{k!}
  \\&= \underbrace{\int_{\mathbb{Z}_p}\cdots\int_{\mathbb{Z}_p}}_{r-\text{times}}    e^{(x_1+\cdots x_r )t}  d\mu_{-1}(x_1)\cdots d\mu_{-1}(x_r)
\end{split}
\end{equation}
Therefore, by (52) and (53) , we obtain the following theorem.
\begin{theorem}
For $p(x) \in \mathbb{P}$, we have
\begin{equation*}
 \langle \underbrace{\int_{\mathbb{Z}_p}\cdots\int_{\mathbb{Z}_p}}_{r-\text{times}} e^{(x_1+\cdots x_r )t}  d\mu_{-1}(x_1)\cdots d\mu_{-1}(x_r)| p(x)\rangle
 =\underbrace{\int_{\mathbb{Z}_p}\cdots\int_{\mathbb{Z}_p}}_{r-\text{times}} p (x_1+\cdots x_r  )  d\mu_{-1}(x_1)\cdots d\mu_{-1}(x_r).
\end{equation*}
Moveover,
\begin{equation*}
\langle\Bigl(\frac{2}{e^t+1}\Bigl)^r|p(x)\rangle=\underbrace{\int_{\mathbb{Z}_p}\cdots\int_{\mathbb{Z}_p}}_{r-\text{times}} p (x_1+\cdots x_r  )  d\mu_{-1}(x_1)\cdots d\mu_{-1}(x_r).
\end{equation*}
\end{theorem}


ACKNOWLEDGEMENTS. This research was supported by Basic Science Research Program through the National Research Foundation of Korea(NRF)
funded by the Ministry of Education, Science and Technology 2012R1A1A2003786.

\par\bigskip

\par

\bigskip\bigskip

\begin{center}\begin{large}

{\sc References}

\end{large}\end{center}

\par

\begin{enumerate}

\item[{[1]}] S. Araci, D. Erdal, J. J. Seo, {\it A study on the fermionic $p$-adic $q$-integral representation on $\Bbb Z_p$  associated with weighted $q$-Bernstein and q-Genocchi polynomials}, Abstr. Appl. Anal. 2011(2011)  Article ID{\bf649248}, 10 pp.

\item[{[2]}] S. Araci, M. Acikg$\ddot{o}$z, H. Jolany, J. J. Seo, {\it A unified generating function of the $q$-Genocchi polynomials with their interpolation functions}, Proc. Jangjeon Math. Soc. {\bf15} (2012), no. 2, 227--233.

\item[{[3]}] L. Carlitz, {\it Some remarks on the multiplication theorems for the Bernoulli and Euler polynomials}, Glas. Mat. Ser. III {\bf16(36)} (1981), no. 1, 3--23.

\item[{[4]}] L. Carlitz, {\it Some polynomials related to the Bernoulli and Euler polynomials}, Utilitas Math. {\bf19} (1981), 81--127.

\item[{[5]}] M. Can, M. Cenkci, V. Kurt, Y. Simsek, {\it On the higher-order $w$-$q$-Genocchi numbers}, Adv. Stud. Contemp. Math. {\bf19} (2009), no. 1, 39--57.

\item[{[6]}] J. Choi, D. S. Kim, T. Kim, Y. H. Kim, {\it Some arithmetic identities on Bernoulli and Euler numbers arising from the $p$-adic integrals on $\Bbb Z\sb p$}, Adv. Stud. Contemp. Math. {\bf22} (2012), no. 2, 239--247

\item[{[7]}] R. Dere, Y. Simsek, {\it Applications of umbral algebra to some special polynomials,} Adv. Stud. Contemp. Math. {\bf22}  (2012), no. 3, 433-438.

\item[{[8]}]  L. C. Jang, {\it A study on the distribution of twisted $q$-Genocchi polynomials}, Adv. Stud. Contemp. Math. {\bf18} (2009), no. 2, 181--189.

\item[{[9]}] D. S. Kim, D.V. Dolgy, T. Kim, S.-H. Rim, {\it Some formulae for the product of two Bernoulli and Euler polynomials}, Abstr. Appl. Anal. 2012(2012), Article ID {\bf784307}, 15 pp.

\item[{[10]}] T. Kim, {\it  Identities involving Frobenius?Euler polynomials arising from non-linear differential equations}, J. Number Theory {\bf132} (2012), no. 12, 2854--2865.

\item[{[11]}] T. Kim,  {\it Symmetry of power sum polynomials and multivariate fermionic $p$-adic invariant integral on $\Bbb Z_p$}, Russ. J. Math. Phys. {\bf16}(2009), no. 1, 93-96 .

\item[{[12]}] S.-H. Rim, J. Jeong, {\it On the modified $q$-Euler numbers of higher order with weight}, Adv. Stud. Contemp. Math. {\bf22} (2012), no. 1, 93-98.

\item[{[13]}]  S.-H. Rim, S.-J. Lee,  {\it Some identities on the twisted $(h,q)$-Genocchi numbers and polynomials associated with $q$-Bernstein polynomials}, Int. J. Math. Math. Sci. 2011(2011), Art. ID {\bf482840}, 8 pp.

\item[{[14]}]  S. Roman, {\it The umbral calculus}, Dover Publ. Inc. New York, 2005.

\item[{[15]}]  Y. Simsek, O. Yurekli, V. Kurt, {\it On interpolation functions of the twisted generalized Frobenius-Euler numbers}, Adv. Stud. Contemp. Math.  {\bf15} (2007), no. 2, 187--194.

\item[{[16]}]  C. S. Ryoo,  {\it Some relations between twisted $q$-Euler numbers and Bernstein polynomials}, Adv. Stud. Contemp. Math. {\bf21} (2011), no. 2, 217--223.

\item[{[17]}] Y. Simsek, O. Yurekli, V. Kurt, {\it On interpolation functions of the twisted generalized Frobenius-Euler numbers}, Adv. Stud. Contemp. Math.(2007), {\bf15}, no. 2, 187-194 .

\item[{[18]}] K. Shiratani,  {\it On Euler numbers}, Mem. Fac. Sci., Kyushu Univ., Ser. A {\bf27}(1973), 1-5 .

\item[{[19]}]   K. Shiratani, S. Yamamoto, {\it On a $p$-adic interpolation function for the Euler numbers and its derivatives},  Mem. Fac. Sci. Kyushu Univ. Ser. A {\bf39} (1985), no. 1, 113--125.

\end{enumerate}

\par

\bigskip\bigskip

\par

\bigskip\bigskip

\par\noindent

\bpn {\small Taekyun {\sc Kim} \mpn Department of Mathematics,
\pn Kwangwoon University, Seoul 139-701, Republic of Korea \pn {\it E-mail:}\ {\sf tkkim@kw.ac.kr} }

\mpn { \bpn {\small Dae San {\sc Kim} \mpn Department of
Mathematics, \pn Sogang University,  Seoul 121-742, Republic of Korea \pn {\it E-mail:}\ {\sf dskim@sogang.ac.kr} }

\mpn { \bpn {\small Sangtae  {\sc Jeong} \mpn Department of
Mathematics, \pn Inha University,  Incheon 860-7114, Republic of Korea \pn {\it E-mail:}\ {\sf stj@inha.ac.kr} }

\bpn {\small Seog-Hoon  {\sc Rim} \mpn  Department of Mathematics Education,\pn Kyungpook National University, Taegu 702-70, Republic of Korea \pn {\it E-mail:}\ {\sf shrim@knu.ac.kr} }

\end{document}